\documentclass{amsproc}

\newtheorem{Theorem}{Theorem}

\newtheorem{theorem}{Theorem}[section]

\newtheorem{proposition}[theorem]{Proposition}
\newtheorem{corollary}[theorem]{Corollary}
\theoremstyle{definition}

\theoremstyle{remark}

\numberwithin{equation}{section}

\newcommand{\R}{{\mathbb R}}
\newcommand{\Z}{{\mathbb Z}}
\newcommand{\FF}{{\mathcal F}}
\newcommand{\GG}{{\mathcal G}}

\newcommand{\N}{{\mathbb N}}

\newcommand{\Fix}{{\rm Fix}}

\newcommand{\DR}{{\rm Diff}_c^r({\mathbb R}^n)_0}
\newcommand{\DS}{{\rm Diff}^r_+(S^1)}
\newcommand{\DSO}{{\rm Diff}^0_+(S^1)}
\newcommand{\DSp}{{\rm Diff}^p(S^1)}

\newcommand{\DB}{{\rm Diff}_c^r(B_i)_0}
\newcommand{\DIp}{{\rm Diff}^p([0,1))}
\newcommand{\DRp}{{\rm Diff}^p({\mathbb R})}

\begin{document}
\title
{New proofs of theorems of Kathryn Mann}


\author[S. Matsumoto]{Shigenori Matsumoto\\ \\{\em To the memory 
of Akio Hattori}}
\address{Department of Mathematics, College of
Science and Technology, Nihon University, 1-8-14 Kanda, Surugadai,
Chiyoda-ku, Tokyo, 101-8308 Japan
}
%\curraddr{Department of Mathematics, College of
%Science and Technology, Nihon University, 1-8-14 Kanda, Surugadai,
%Chiyoda-ku, Tokyo, 101-8308 Japan}
\email{matsumo@math.cst.nihon-u.ac.jp
}
\thanks{The author is partially supported by Grant-in-Aid for
Scientific Research (C) No.\ 25400096.}
%    General info


\date{\today }
\begin{abstract}
We give a shorter proof of the following theorem of
Kathryn Mann \cite{M}: the identity component of the
group of the compactly supported $C^r$ diffeomorphisms
of $\R^n$ 
cannot admit a nontrivial $C^p$-action on $S^1$,
provided $n\geq2$,  $r\neq n+1$ and $p\geq2$.
We also give a new proof of another theorem of Mann \cite{M}: any nontrivial
 homomorphism from the group
of the orientation preserving $C^r$ diffeomorphisms of the circle 
to the group of $C^p$ diffeomorphisms of the circle
is the conjugation of the standard inclusion 
by a $C^p$ diffeomorphism, if $r\geq p$, $r\neq2$ and $p\neq1$.
\end{abstract}

\maketitle

\section{Introduction}

\'E. Ghys \cite{G} asked if the group of diffeomorphisms of a manifold
admits a nontrivial action on a lower dimensional manifold.
A break through towards this problem was obtained by Kathryn Mann \cite{M}
in the case where the target manifold is one dimensional.
Let us denote by $\DR$ the identitiy component of the
group of the compactly supported $C^r$ diffeomorphisms
of $\R^n$,  $r=0,1,\ldots,\infty$.
She showed the following theorem.

\begin{Theorem} \label{T}
Assume $n\geq2$, $r\neq n+1$ and $p\geq2$. Then any (abstract)
homomorphism from $\DR$ to $\DSp$ is trivial.
\end{Theorem}

The condition $r\neq n+1$ is for the simplicity of the source
group.
The condition $p\geq2$ is necessary since the proof
is built upon a theorem of Kopell
and Szekeres.
Notice that by the fragmentation lemma, the same statement holds true if
we replace $\R^n$ by any $n$-dimensional manifold,
compact or not.
One aim of this notes is to give a short proof of 
the above theorem. We also show the following result.

\begin{Theorem}\label{T2}
Assume $n\geq2$,  $r\neq n+1$ and $p\geq2$. Then any 
homomorphism from $\DR$ to $\DRp$ is trivial.
\end{Theorem}

This is a generalization of a result of \cite{M}
for the target group $\DIp$.

\medskip

Next we consider the case where the source manifold
is one dimensional. We provide a shorter proof of the following theorem,
also contained in \cite{M}. Denote by $\DS$ the group of
the orientation preserving $C^r$ diffeomorphisms  of $S^1$.

\begin{Theorem} \label{T'}
Assume $r\geq p$, $r\neq2$ and $p\neq1$. Then any nontrivial homomorphism from $\DS$
to $\DSp$ is
the conjugation of the standard inclusion by a $C^p$ diffeomorphism.
\end{Theorem}

In the above theorem, the case where $p=0$ is new. We also have the following.

\begin{Theorem} \label{T''}
Assume $p\neq1$. Then any nontrivial homomorphism from ${\rm PSL}(2,\R)$ to ${\rm Diff}^p(S^1)$
is the conjugation of the standard inclusion by a $C^p$ diffeomorphism.
\end{Theorem}

As for Ghys's question for target manifolds of dimension $>1$, a satisfactory
answer is obtained by S. Hurtado \cite{H}. Some part of
his argument is an induction on the dimension of the target
manifold. It is
based upon Theorems \ref{T} and \ref{T'}.

In \cite{M}, Theorems \ref{T} and \ref{T'} are shown using the following result.

\begin{theorem} \label{ht}
Assume $r\geq 3$, $p\geq2$ and $r\geq p$. Any nontrivial homomorphism $\Phi$
from ${\rm Diff}^r_c((0,1))$
to ${\rm Diff}^p([0,1))$ without interior
global fixed point of the $\Phi$-action
is the conjugation of the standard inclusion
by a $C^p$ diffeomorphism of $(0,1)$.
\end{theorem}

Our  proofs of Theorems \ref{T} and \ref{T'} do not
use Theorem \ref{ht}. On the other
hand, we would like to stress that
Theorem \ref{ht} is more involved, and cannot be shown 
by the techniques of the present paper.

\section{Theorem of Kopell and Szekeres}

Our main tool for the proof of
Theorems \ref{T} and \ref{T2} is the following theorem due to
Kopell \cite{K} and Szekeres. (See 4.1.11 in \cite{N2}.)
This forces us to assume $p\geq2$ in these theorems. 

\begin{theorem} \label{t2}
 Let $p\geq2$. Assume that
$\psi\in{\rm Diff}^p([0,1))$ admits no interior fixed point.
Then there is a unique $C^1$ flow $\{\psi^t\}$ on $[0,1)$
such that $\psi=\psi^1$.
Moreover any element $\phi$ of the centralizer $C(\psi)$ of $\psi$ in
${\rm Diff}^p([0,1))$ can be written as $\phi=\psi^t$ for some $t\in\R$.
\end{theorem}

\begin{corollary}\label{c}
Let $p\geq2$ and $\psi\in{\rm Diff}^p([0,1))$. Then $C(\psi)=C(\psi^2)$.
\end{corollary}

{\sc Proof}. Choose any element $g\in C(\psi^2)$. 
Let $J$ be
the closure of any component of $[0,1)\setminus\Fix(\psi)$.
(Notice that $\Fix(\psi^2)=\Fix(\psi)$.)
Then by Theorem \ref{t2}, $g$ commutes with $\psi$ on $J$.
Since $J$ is arbitrary, $g$ commutes with $\psi$
everywhere.
\qed

\medskip
We also have the following result, whose proof is the same
as above.

\begin{corollary}\label{c2}
Assume $p\geq2$ and $\psi\in{\rm Diff}^p_+(\R)$  
 admits
fixed points. Then the centralizers in
${\rm Diff}^p_+(\R)$ satisfy $C(\psi)=C(\psi^2)$.
\end{corollary}

\section{Commuting subgroups of $\DSO$}

Another basic result needed for the proof is the following.

\begin{proposition} \label{p1}
Let $G_1$ and $G_2$ be simple nonabelian
subgroups of $\DSO$.
Assume that $G_2$ is conjugate to $G_1$ in
$\DSO$ and that any element of $G_1$ commutes with any element of $G_2$.
Then there is a global fixed point of $G_1$:
$\Fix(G_1)\neq\emptyset$.
\end{proposition}

{\sc Proof}. 
First of all, let us show that there is an element 
$g\in G_1\setminus\{{\rm id}\}$
such that $\Fix(g)$ is nonempty.
Assume the contrary. Consider the group $\tilde G_1$ formed
by any lift of any element of $G_1$ to the universal covering space
$\R\to
S^1$.
The canonical projection
$\pi:\tilde G_1\to G_1$ is a group homomorphism.
Now $\tilde G_1$ acts freely on $\R$. 
A theorem of H\"older asserts that $\tilde G_1$ is abelian.
See for example \cite{N1}. Therefore $G_1=\pi(\tilde G_1)$
would be abelian, contrary to the assumption of the proposition.
 
Let $X_2\subset S^1$ be a minimal set of $G_2$.
The set $X_2$ is either a finite set, a Cantor set or
the whole of $S^1$. If $X_2$ is a singleton, then
$G_2$ admits a fixed point. Since $G_1$ is conjugate to
$G_2$, we have $\Fix(G_1)\neq\emptyset$, as is required.
If $X_2$ is a finite set which is not a singleton, we
get a nontrivial homomorphism from $G_2$ to a
finite abelian group, contrary to the assumption.
In the remaining case, it is well known, easy to show,
that the minimal set is unique. That is, $X_2$
is contained in any nonempty $G_2$ invariant closed subset.

Let $F_1$ be the subset of $G_1$ formed by the elements with nonempty
fixed point set. 
Since $G_1$ and $G_2$ commutes, the fixed point set
$\Fix(g)$ of any element $g\in F_1$ is $G_2$ invariant.
 Then we have:
\begin{equation}\label{ee}
\mbox{ $X_2\subset\Fix(g)$ for any $g\in F_1$}.
\end{equation}
This shows that $F_1$ is 
in fact a subgroup. By the very definition,
$F_1$ is normal. Since $G_1$ is simple and $F_1$ is nontrivial, $F_1=G_1$.
Finally again by (\ref{ee}), $\Fix(G_1)\neq\emptyset$.
\qed

\section{Proof of Theorem \ref{T}}

Assume $n\geq2$, $r\neq n+1$ and $p\geq2$.
 Let $\Phi:\DR\to\DSp$ be a nontrivial homomorphism.
Our purpose is to deduce a contradiction.
Since $\DR$ is simple by the assumption $r\neq n+1$, 
the map $\Phi$ is injective and its image is contained
in the group of orientation preserving diffeomorphisms.

Let $B_1$ and $B_2$ be disjoint open balls of radius 2 in $\R^n$.
The group $\GG_i=\DB$ is nonabelian and simple.
Clearly $\GG_2$ is conjugate to $\GG_1$
in $\DR$.

Therefore $\Phi(\GG_i)$ satisfies all the conditions of
Proposition \ref{p1}. Thus $\Phi(\GG_1)$ has a fixed point,
and one can identify $\Phi(\GG_1)\subset \DIp$.
In view of Corollary \ref{c} and the injectivity of $\Phi$, 
it is sufficent to construct an element
$g\in\GG_1$ such that $C(g)\neq C(g^2)$.
Let $B_1'$ be the concentric ball in $B_1$ of radius $1$.
Any element $g\in\GG_1$ which is an involution on $B_1'$
will do. 

\section{Proof of Theorem \ref{T2}}

Let $n\geq2$, $r\neq n+1$
and $p\geq2$. Assume there is a nontrivial homomorphism $\Phi:\DR\to\DRp$.
By the simplicity of $\DR$, $\Phi$ is injective, with its image
contained
in the group of orientation preserving diffeomorphisms.
In view of Corollary \ref{c2} and the 
last step of  the previous section, it suffices to show that any element
of the image of $\Phi$ has nonempty fixed point set.
The rest of this section is devoted to its proof.

Assume for contradiction that there is an element $g'\in\DR$
such that $\Fix(\Phi(g'))=\emptyset$.
Choose open balls $B_i$ ($i=1,2$) in $\R^n$ as in the previous
section. Again let $\GG_i={\rm Diff}^r_c(B_i)_0$.
There is a conjugate $g$ of $g'$ in $\GG_2$.
Notice that  $\Fix(\Phi(g))=\emptyset$. 
Then $\Phi(\GG_2)$
has a cross section $I$ in $\R$, that is, $I$ is a compact interval
such that any $\Phi(\GG_2)$
orbit  hits $I$.
Now we follow the proof of Proposition 6.1 in \cite{DKNP}, to show that 
there is a unique minimal set $X_2$ for $\Phi(\GG_2)$.
Moreover we shall show that there is
a nonempty $\Phi(\GG_2)$ invariant closed subset $X_2$ in $\R$ which
has the property that 
any nonempty $\Phi(\GG_2))$ invariant closed  subset
contains $X_2$.

The proof goes as follows. 
Let $F$ be the family 
of nonempty $\Phi(\GG_2)$ invariant closed  subsets of $\R$,
and  $F_I$ the family of nonempty closed subsets $Y$ in $I$ such that
 $\Phi(\GG_2)(Y)\cap I=Y$,
where we denote
$$\Phi(\GG_2)(Y)=\bigcup_{g\in\GG_2}\Phi(g)(Y).$$
Define a map $\phi:F\to F_I$ by $\phi(X)=X\cap I$,
and $\psi:F_I\to F$ by $\psi(Y)=\Phi(\GG_2)(Y)$.
They satisfiy $\psi\circ\phi=\phi\circ\psi={\rm id}$.

Let $\{Y_\alpha\}$ be a totally ordered set in $F_I$.
Then the intersection $\cap_\alpha Y_\alpha$ is nonempty.
Let us show that it belongs to $F_I$, namely,

\begin{equation}\label{e2}
\Phi(\GG_2)(\cap_\alpha Y_\alpha)\cap I=\cap_\alpha Y_\alpha.
\end{equation}
For the inclusion $\subset$, we have
$$
\Phi(\GG_2)(\cap_\alpha Y_\alpha)\cap I\subset 
(\cap_\alpha \Phi(\GG_2)(Y_\alpha))\cap I=
\cap_\alpha (\Phi(\GG_2)(Y_\alpha)\cap I)=\cap_\alpha Y_\alpha.
$$
For the other inclusion, notice that
$$\cap_\alpha Y_\alpha\subset \Phi(\GG_2)(\cap_\alpha Y_\alpha)
\mbox{ and }\cap_\alpha Y_\alpha\subset I.
$$

Therefore by Zorn's lemma, there is a minimal element $Y_2$ in $F_I$. 
The set $Y_2$ is not finite. In fact, if it is finite,
the set $X_2=\psi(Y_2)$ in $F$ is discrete, and 
there would be a nontrivial homomorphism from $\Phi(\GG_2)$
to $\Z$, contrary to the fact that $\GG_2$, (and hence $\Phi(\GG_2)$)
is simple.

Now the correspondence $\phi$ and $\psi$ preserve the inclusion.
This shows that  
there is no nonempty $\Phi(\GG_2)$
invariant closed proper subset of $X_2=\psi(Y_2)$.
In other words, any $\Phi(\GG_2)$ orbit 
contained in
 $X_2$ is dense in $X_2$.
Therefore $X_2$ is either $\R$ itself or a locally
Cantor set. 
In the former case, any nonempty $\Phi(\GG_2)$ invariant
closed subset must be $\R$ itself. 
 
Let us show that in the latter case,
$X_2$ satisfies the desired property: $X_2$ is contained in
any nonempty $\Phi(\GG_2)$ invariant closed subset.
For this, we only need to show that
the  $\Phi(\GG_2)$ orbit of any point $x$ in $\R\setminus X_2$ 
accumulates to a point in $X_2$. Let $(a,b)$ be the connected component
of $\R\setminus X_2$ that contains $x$. Then
there is a sequence $g_k\in\GG_2$ ($k\in\N$) 
such that $\Phi(g_k)(a)$ accumulates
to $a$ and that $\Phi(g_k)(a)$'s are mutually distinct. Then the intervals
$\Phi(g_k)((a,b))$ are mutually disjoint, and consequently
$\Phi(g_k)(x)$ converges to $a$.
This concludes the proof that $X_2$ is contained in
any nonempty $\Phi(\GG_2))$ invariant closed subset.

\medskip

Let $\FF_1$ be the subset of $\GG_1$ formed by the elements $g$
such that $\Fix(\Phi(g))\neq\emptyset$. 
Again by a theorem of H\"older, $\FF_1$ contains a nontrivial element. 
Now we have
\begin{equation}\label{e1}
\mbox{ $X_2\subset\Fix(\Phi(g))$ for any $g\in \FF_1$}.
\end{equation}
This shows that $\FF_1$ is 
a subgroup, normal by the definition.
Since $\GG_1$ is simple, $\FF_1=\GG_1$.
Finally again by (\ref{e1}), we have $\Fix(\Phi(\GG_1))\neq\emptyset$.
This contradicts the fact that $\Phi(\GG_2)$, being
conjugate to $\Phi(\GG_1)$, must also
have a free element. The contradiction shows that
 any element
of the image $\Phi$ has nonempty fixed point set.

\section{Proof of Theorem \ref{T'}}

We first prove Theorem \ref{T'} for $p=0$. 
Assume  $r\neq2$ and let
$\Phi$ be a nontrivial homomorphism from
$F=\DS$ to ${\rm Diff}^0(S^1)$.
Since $F$ is simple, $\Phi$ is injective and the image of $\Phi$
is contained in $\DSO$.
For any $x\in S^1$, denote by $F_x$ the isotropy subgroup
at $x$. 

\begin{proposition}\label{p3}
For any $x\in S^1$, the fixed point set $\Fix(\Phi(F_x))$ is a singleton.
\end{proposition}

The proof uses Theorem 5.2 in \cite{M'}, which states as follows.

\begin{theorem} \label{t3}
Any nontrivial homomorphism from ${\rm PSL}(2,\R)$ to ${\rm Diff}^0_+(S^1)$
is the conjugation of the standard inclusion by a homeomorphism $h$.
\end{theorem}

{\sc Proof of Proposition \ref{p3}}. 
Let $G_x={\rm Diff}^r_c(S^1\setminus\{x\})$.
First we shall show $\Fix(\Phi(G_x))\neq\emptyset$.
Let $\{U_n\}_{n\in\N}$ be an increasing sequence of precompact
open intervals in $S^1\setminus\{x\}$ such that
$\cup_nU_n=S^1\setminus\{x\}$ and let $V_n$ be an open
interval in $S^1\setminus\{x\}$ disjoint from $U_n$.
Then $\Phi({\rm Diff}^r_c(U_n))$ and $\Phi({\rm Diff}^r_c(V_n))$
satisfy the conditions of Proposition \ref{p1}. 
Therefore $\Fix(\Phi({\rm Diff}^r_c(U_n)))$ is nonempty.
Since $\Phi(G_x)=\cup_n\Phi({\rm Diff}^r_c(U_n))$, we have
$\Fix(\Phi(G_x))=\cap_n\Fix(\Phi({\rm Diff}^r_c(U_n)))$.
But the RHS is a decreasing intersection of nonempty compact subsets.
Therefore we get $\Fix(\Phi(G_x))\neq\emptyset$.

Next let us show that $\Fix(\Phi(G_x))$ is a singleton.
Assume there are two distinct points $\xi_0$ and $\xi_1$ in
$\Fix(\Phi(G_x))$.
By Theorem \ref{t3}, there is a rotation $R\in{\rm PSL}(2,\R)$
such that $\Phi(R)(\xi_0)=\xi_1$. Then $\xi_1$ is left fixed both
by $\Phi(G_x)$ and $\Phi(RG_xR^{-1})=\Phi(G_{R(x)})$,
and hence by $\Phi(F)$, since $G_x$ and $G_{R(x)}$ generate $F$.
Especially $\Phi({\rm PSL}(2,\R))$ admits a global
fixed point, contradicting Theorem \ref{t3}.

Finally since $G_x$ is a normal subgroup of $F_x$ and since 
$\Fix(G_x)$ is a singleton, we have
$\Fix(F_x)=\Fix(G_x)$.
\qed

\bigskip
Define a map $h':S^1\to S^1$ by $\{h'(x)\}=\Fix(F_x)$.
Let us show that $h'$ coincides with a homeomorphism $h$
in Theorem \ref{t3}. In fact, for any $x\in S^1$, consider
a parabolic element $g\in{\rm PSL}(2,\R)$ such that $g(x)=x$.
Then, by Theorem \ref{t3}, $h(x)$ is the unique fixed point of $\Phi(g)$.
Therefore $h(x)=h'(x)$.

Now for any $f\in F$, we have
\begin{eqnarray*}
\{h(f(x))\}&=&\Fix(\Phi(F_{f(x)}))=\Fix(\Phi(fF_xf^{-1}))
=\Fix(\Phi(f)\Phi(F_x)\Phi(f)^{-1})\\
&=&\Phi(f)(\Fix(\Phi(F_x))
=\Phi(f)\{h(x)\}.
\end{eqnarray*}
This shows that $h\circ f=\Phi(f)\circ h$ for any $f\in F$.
Namely $\Phi$ is the conjugation by a homeomorphism $h$.
This completes the proof for $p=0$.

\medskip
Let us consider the case $p\geq 2$.
We follow an argument in \cite{T},
and show that the homeomorphism $h$ established above
is in fact a $C^p$ diffeomorphism.
Denote $H={\rm Diff}^p_+(S^1)$. 
First of all, since $h$ is locally monotone and thus
of bounded variation, it is
differentiable Lebesgue almost everywhere. But we have
$h \circ R_t=f_t\circ h$ ($\forall t\in S^1$), where
$R_t$ is the rotation by $t$ and $f_t\in H$.
This shows that $h$ is differentiable everywhere, with nonvanishing
derivative.

Let us show that for any point $x\in S^1$, $h$ is a $C^p$ diffeomorphism
on some neighbourhood of $x$. Choose $f\in F$ such that $f(x)=x$ and
$f'(x)=\lambda\in(0,1)$. 
Then $h\circ f\circ h^{-1}\in H$ leaves
$h(x)$ fixed and the derivative there is also $\lambda$.
Notice that $f$ also belongs to $H$, since $r\geq p$.
By the Sternberg linearization theorem (Theorem 2 of \cite{S}),
there is a $C^p$ diffeomorphism $k_1$ (resp.\ $k_2$) from a neighbourhood
of $x$ (resp. $h(x)$) to a neighbourhood of $0$ in $\R$ such that
$k_1\circ f=L\circ k_1$ (resp. $k_2\circ (h\circ f\circ h^{-1})=L\circ
k_2$),
where $L$ is a linear map of $\R$ with slope $\lambda$.
(The assumption $p\neq1$ is necessary for Theorem 2 of \cite{S}.)
Then the composite $k_2\circ h\circ k_1^{-1}$ is a local homeomorphism
of $(\R,0)$ commuting with $L$ and differentiable at $0$.
It is easy to prove that $k_2\circ h\circ k_1^{-1}$ is a linear map.
(The linear contraction $L\times L$ of $\R^2$ leaves the
graph of $k_2\circ h\circ k_1^{-1}$ invariant.)
Then $h$ is a $C^p$ diffeomorphism in a neighbourhood of $x$, as is required.

\bigskip
The proof of Theorem \ref{T''} is given by the same argument as above
starting from the homomorphism $h$ in Theorem \ref{t3}.

\end{document}